\documentclass[12pt, fleqn]{article}
\usepackage{amssymb,amsmath,amsfonts, ,amsthm, graphicx, color}
\usepackage{float}
\usepackage{rotating}
\usepackage{lipsum}
\newcommand{\NI}{\noindent}
\newtheorem{theorem}{\NI{\bf Theorem}}[section]

\newtheorem{cor}{\NI\bf Corollary}[section]
\newtheorem{defn}{\NI\bf Definition}[section]

\newcommand{\bp}{\begin{proof}}
\newcommand{\ep}{\end{proof}}
\newcommand{\bt}{\begin{theorem}}
\newcommand{\et}{\end{theorem}}
\newcommand{\bc}{\begin{cor}}
\newcommand{\ec}{\end{cor}}
\newcommand{\newsection}[1]{\setcounter{equation}{0}
\setcounter{theorem}{0}
\section{#1}}

\newcommand{\bea}{\begin{eqnarray}}
\newcommand{\eea}{\end{eqnarray}}
\newcommand{\ben}{\begin{eqnarray*}}
\newcommand{\een}{\end{eqnarray*}}

\def \qed {\hfill \vrule height 6pt width 6pt depth 0pt}
\newcommand{\be}{\begin{equation}}
\newcommand{\ee}{\end{equation}}

\title {\bf Testing Exponentiality Against a Trend Change in Mean Time to Failure in Age Replacement}
\author{
Muhyiddin Izadi\thanks {Corresponding author} $^{\ **}$, \ \  Sirous Fathimanesh$^{***}$
\\  \\ $^{**}$Department of Statistics, Razi
University, Kermanshah, Iran \\
$^{***}$Department of Statistics, University of Kurdistan, Sanandaj, Iran
}
\date{}
\begin {document}
\maketitle
\begin{abstract}
Mean time to failure in age replacement evaluates the performance and effectiveness of the age replacement policy. 
In this paper, we propose a test for exponentiality against a trend change in mean time to failure in age replacement. We derive the asymptotic distribution of the test statistics under the null hypothesis to approximate the critical values. We conduct a simulation study to investigate the performance of the proposed test and compare it with some well known tests in the literature.

\NI {{\bf Keywords} : BFR distribution,  Durbin's approximation, Mean time to failure,  Non-monotonic aging class, NWBUE distribution, Total time on test transform.}
\end{abstract}

\newsection{Introduction}

The failure of a component or a system during   is usually costly or dangerous. It is a common practice to employ a maintenance policy to prevent the item from the failure during operation. The most common and popular maintenance policy is the age replacement policy in which  an item is replaced by a new one upon failure or at a known age $t$, whichever comes first.   Let $F$ be the lifetime distribution of a new item and $X_{[t]}$ denote the time to the first in-service failure of an item under the age replacement policy with the age replacement time $t$.  Then the reliability function of $X_{[t]}$ (denoted by $R_t$ )  is given by 
\begin{eqnarray*}
R_t(x)= [\overline{F}(t)]^n\overline{F}(x-nt),    \ \ \  nt \le x< (n+1)t, \ \ \ n=0, 1, \ldots,
\end{eqnarray*}
 where $\overline{F}=1-F$ (Barlow and Proschan, 1965).   To evaluate the performance and   effectiveness of the age replacement policy, Barlow and Proschan (1965) introduced the mean of $X_{[t]}$, which is called the mean time to failure (MTTF) in the age replacement policy.    Let us denote $M_F(t)=E(X_{[t]})$. It is known that
 \begin{eqnarray}\label{MTTF}
 M_F(t)=\frac{\int_{0}^{t}\overline{F}(x)dx}{F(t)},  \ \ \ t>0
 \end{eqnarray}
 (Kayid et al. 2013; Izadi et al., 2018). 
  
 Study  the  behaviour  of  $M_F(t)$  with  respect  to  $t$  might lead  us  to  realize  the  optimal
replacement time which makes the MTTF of practical importance.   The distribution $F$ is said to be decreasing mean time to failure (DMTTF) in age replacement if   $M_F(t)$ is decreasing in $t \in (0, \infty)$ which means a kind of ``deterioration". The dual class of increasing  mean time to failure (IMTTF) in age replacement is defined by changing the sense of the monotonicity  and means  ``non-deterioration" or improvement in some senses.   The ageing classes DMTTF and IMTTF have been  studied in the literature by many
researchers. Klefsj\"{o} (1982) investigated the relationship between DMTTF and IMTTF and some well known aging classes of distributions. 
 Knopik (2005) showed that the DMTTF and IMTTF classes are closed
under formation of parallel and series systems.  Knopik (2006) further studied the relationship between the DMTTF (IMTTF) class and the increasing (decreasing)  failure rate in average (IFRA (DFRA)) class of distributions and showed that the DMTTF is closed under weak convergence of distribution and convolution. Li and Xu (2008) introduced the NBUR$_{rh}$ class of
life distributions which is equivalent to the DMTTF class. They studied various properties of the DMTTF class and provided a test for exponentiality against monotone MTTF. Asha and Nair (2010) studied some properties of the quantile based MTTF function by examining its relationship  with hazard (reversed hazard) rate and mean (reversed mean) residual life functions. They also defined a new stochastic ordering of life distributions based on MTTF and studied its relationship with some known orderings useful in reliability analysis.   Kayid et al. (2013) studied the preservation properties of the MTTF order under  monotonic transformations, mixture,
and weighted distributions. The problem of  testing exponentiality against the DMTTF property has been considered in Kayid et al. (2013) and Kattumannil and Anisha (2016).



The ageing patterns in the above classes are monotone. But, in practical situations, it is often seen that the ageing
pattern is non-monotonic. In order to model such situations, various non-monotonic ageing classes have been defined in the literature, see for example, Glaser (1980), Rajarshi and Rajarshi (1988), Klefsj\"{o} (1989), Deshpande and Suresh (1990), Mitra and Basu (1994),  Belzunce et al. (2007). 

Izadi et al.  (2018) proposed   the following two  non-parametric  classes  of  distributions with non-monotonic MTTF function.
\begin{defn}
A life distribution $F$  is said to be an increasing  then decreasing mean time to failure (IDMTTF) (decreasing then  increasing  mean time to failure (DIMTTF)) distribution if there exists a turning point $\tau\ge0$ such that $M_F(t)$ is increasing (decreasing) on $(0, \tau]$ and decreasing (increasing) on $(\tau, \infty)$. 
\end{defn}

The  IDMTTF class of distributions models  a situation that the effect of age replacement is initially  beneficial and  then adverse and the DIMTTF class of  distributions can be used to  model the cases when the effect of age replacement is initially adverse and then beneficial.  In the IDMTTF class of distributions, $M_F(t)$ is maximum at the change point, so it may be taken as a possible optimal age replacement. Thus, the change point of $M_F(t)$ is important in the IDMTTF distributions which makes this class of distributions of great interest in connection with the age replacement optimization. One of the most used criteria to determine the optimal replacement time is minimizing the expected cost rate (cf. Nakagawa , 2005).  Izadi et al. (2018) studied the relationship between the MTTF function and the expected cost rate function. They investigated how $\tau$ approximates the optimal replacement time which minimizes the expected cos rate function. Izadi et al. (2018) also studied the implications between the IDMTTF and DIMTTF classes
of distributions and some existing classes of non-monotonic aging classes. 

Let $F$ be a lifetime distribution with finite mean $\mu_F$, density function $f$ and failure rate function $r_F(t)=\frac{f(t)}{\overline{F}(t)}$.  It is said that $F$ is bathtub (upside-down) failure rate (BFR (UBFR)), if there exists $t_0\ge0$ such that $r_F(t)$ is decreasing (increasing) in $0\le t< t_0$ and increasing (decreasing) in $t\ge t_0$; $F$ is NWBUE (NBWUE) if there exists $t^*\ge0$ such that $\mu_F \le (\ge) \frac{\int_{t}^{\infty}\overline{F}(x)}{\overline{F}(t)}$ for $t\le t^*$ and $\mu_F \ge (\le) \frac{\int_{t}^{\infty}\overline{F}(x)}{\overline{F}(t)}$ for $t> t^*$.   

Izadi et al. (2018) showed that
\begin{eqnarray*}
\mbox{BFR}  \Longrightarrow IDMTTF \Longrightarrow NWBUE \ \mbox{ and } \mbox{UBFR}  \Longrightarrow DIMTTF \Longrightarrow NBWUE.
\end{eqnarray*} 

A reasonable starting point in reliability analysis is to determine the ageing class of the underlying distribution $F$.  In view of this consideration, the statistical problem of testing whether the lifetime distribution $F$   belongs to a specific monotonic or non-monotonic ageing class has been received considerable attention in the literature; see for instance Guess et al. (1986), Klefsj\"{o} (1989), Hawkins et al. (1992), Hawkins and Kochar (1997), Lai (1994), Na et al. (2005),  Lai and Xie (2006), Anis (2014) and Anis and Ghosh (2015) among others. The problem of testing whether $F$ is DMTTF has also been recently studied by Li and Xu (2008), Kayid et al. (2013) and Kattumannil and Anisha (2016).  

Let ${\mathcal{ E}}$ denote the exponential family of distributions, that is, $\mathcal{E}=\{ F| F(x)=\lambda e^{-\lambda x}, \lambda>0\}$ and let 
\begin{eqnarray*}
\mathcal{F}_{ID}&=&\left\{F| F \mbox{ is an IDMTTF distribution and not exponential}\right\}
\end{eqnarray*}
and
\begin{eqnarray*}
\mathcal{F}_{DI}&=&\left\{F| F \mbox{ is an DIMTTF distribution and not exponential}\right\}.
\end{eqnarray*}

 In this paper we are interested the problem of testing
\begin{eqnarray*}
  H_0&:& F \in {\mathcal{ E}}
  \end{eqnarray*}
  against
  \begin{eqnarray}\label{IDMRLH}
  H_1&:& F\in\mathcal{F}_{ID}
\end{eqnarray}
based on the random sample $X_1, X_2, \ldots, X_n$  from distribution $F$.  When  the  dual  model  is  considered, we  test $H_0$ against
\begin{eqnarray}\label{DIMRLH}
H'_1: F\in\mathcal{F}_{DI}.
\end{eqnarray}


The organization of the paper is as follows. In the next section, we propose a test for the considered problem of hypothesis testing  and obtain the asymptotic null distribution of the test statistic. The performance of the proposed test is investigated by simulation study in Section 3. 

\section{Test for trend change}
Let $F$ be a lifetime distribution with finite mean $\mu_F$, density function $f$ and reliability function  $\overline{F}=1-F$.  Taking derivative of $M_F(t)$ with respect to $t$, it can be shown that $F$ is IDMTTF  with the change point $\tau\ge0$, if and only if
\begin{eqnarray*}\label{nq1}
\omega(t)=F(t)\overline F(t)-f(t)\int_{0}^{t}\overline F(x)dx
\left\{%
\begin{array}{ll}
    \ge 0, & t<\tau \\
    \le 0, & t\ge\tau. \\
\end{array}%
\right.
\end{eqnarray*}
Let us define
\begin{eqnarray*}
  \gamma(F,t)&=&
  \int_{0}^{t}\omega(x)dx-\int_{t}^{\infty}\omega(x)dx, \ \ 0\le
  t<\infty
    \end{eqnarray*}
and
\begin{eqnarray*}
 \gamma(F)=\sup\{ \gamma(F,t); 0\le t<\infty\}.
\end{eqnarray*}

Under the null hypothesis $H_0: F$ is an exponential distribution, $\omega(t)=0$ for all $t\ge0$ which implies $\gamma(F)=0$. Taking derivative of $\gamma(F,t)$, it is easy to see that under the
alternative hypothesis $H_1: F$ is IDMTTF and not exponential, $\gamma(F,t)$ is  increasing  for $t<\tau$ and decreasing for $t\ge \tau$ and hence $\gamma(F)=\gamma(F,\tau)>0$. Thus we can use $\gamma(F)$ as a
measure of departure from $H_0$ in favor of $H_1$. It can easily be seen that
    \begin{eqnarray*}
      \gamma(F,t)&=&2(1+\overline F(t))\int_{0}^{t}\overline F(x)dx-2\int_{0}^{\infty}\overline F^2(x)dx+4\int_{t}^{\infty}\overline
      F^2(x)dx-\mu_F\\
      &=& 2(1+\overline F(t))\int_{0}^{\infty}\overline F(x)dx-2\int_{0}^{\infty}\overline F^2(x)dx \\
      &&+\int_{t}^{\infty}\overline F(x) [4\overline F(x)-2(1+\overline F(t))]dx-\mu_F
    \end{eqnarray*}
Now, let  $X_1, \ldots, X_n$ be a  random sample from  $F$ and  $X_{(1)},  \ldots, X_{(n)}$ denote the corresponding order statistics, $F_n$ the empirical distribution function and  $\overline X$ the sample
mean. Set $X_{(0)}\equiv0$ and $X_{(n+1)}\equiv\infty$. 
The test statistic for testing $H_0$ against $H_1$  can be estimate of $\gamma(F)$  by replacing $F$ by $F_n$ as follows: 
\begin{eqnarray*}
  \gamma(F_n)&=& \sup\{ \gamma(F_n,t); 0\le t<\infty\}\\
  &=& \max_{0\le k\le n} \sup\{ \gamma(F_n,t); X_{(k)}\le t<X_{(k+1)}\}.
\end{eqnarray*}
For $X_{(k)}\le t<X_{(k+1)}$, we have
\begin{eqnarray*}
\gamma(F_n,t)&=&2(1+\overline F_n(t))\int_{0}^{\infty}\overline F_n(x)dx-2\int_{0}^{\infty}\overline F_n^2(x)dx-\mu_{F_n}\\
     && +\int_{X_{(k+1)}}^{\infty}\overline F_n(x) [4\overline F_n(x)-2(1+\overline F_n(t))]dx\\&&+\int_{t}^{X_{(k+1)}}\overline F_n(x) [4\overline F_n(x)-2(1+\overline F_n(t))]dx\\
      &=&2(2-\frac{k}{n})\sum_{i=0}^{n-1}(1-\frac in) D_i- \sum_{i=0}^{n-1}(1-\frac in)^2 D_i-\sum_{i=0}^{n-1}(1-\frac in) D_i\\
      && + \sum_{i=k+1}^{n-1} (1-\frac in)[4(1-\frac in)-2(2-\frac kn)]D_i - 2(X_{(k+1)}-t)\frac kn(1-\frac kn)
\end{eqnarray*}
where $D_i=X_{(i+1)}-X_{(i)}, \ \ i=0,1,\ldots, n-1$. Thus, 
 \begin{eqnarray*}
  \sup_{  X_{(k)}\le t<X_{(k+1)}}\gamma(F_n,t)&=&2(2-\frac{k}{n})\sum_{i=0}^{n-1}(1-\frac in) D_i- \sum_{i=0}^{n-1}(1-\frac in)^2 D_i -\sum_{i=0}^{n-1}(1-\frac in) D_i\\
 && + \sum_{i=k+1}^{n-1} (1-\frac in)[4(1-\frac in)-2(2-\frac kn)]D_i 
 \end{eqnarray*}
and hence
\begin{eqnarray*}
  \gamma(F_n)=\max_{0\le k\le
  n}\zeta_{n,k}
\end{eqnarray*}
where 
\begin{eqnarray*}
\zeta_{n,k}=\sum_{i=0}^{n-1}A_nD_{i}+\sum_{i=k+1}^{n-1}B_nD_i,
\end{eqnarray*}
\begin{eqnarray*}
A_n=-2(1-\frac in)^2-(1-\frac in) + 2(2-\frac kn)(1-\frac in), \ B_n=4(1-\frac in)^2-2(2-\frac kn)(1-\frac in).
\end{eqnarray*}
 The IDMTTF class of distributions is closed under the scale transformation, that is, if the random variable $X$ is IDMTTF, then $aX$ is also IDMTTF for every $a>0$.   Since $ \gamma(F_n)$ is not scale invariant, we use the test statistic 
 \begin{eqnarray*}
  \gamma^*(F_n)=\frac{\sqrt{n}\gamma(F_n)}{\overline X}.
\end{eqnarray*}
which makes the test-scale invariant. 

In the next theorem, we obtain the asymptotic distribution of
$\gamma^*(F_n)$ under the null hypothesis.
\begin{theorem}\label{Th1}
  Let $Z(u)$ be a mean zero Gaussian process with covariance
\begin{eqnarray}\label{Covariance}
\rho(s,t)=\frac23(s^3-t^3)+\frac13, \ \ \ \ \ 0\le s\le t\le1.
\end{eqnarray}
  Under $H_0$,
  \begin{eqnarray*}
    \gamma^*(F_n)\stackrel D\rightarrow Z=\sup\{Z(u); 0\le u\le1\},
  \end{eqnarray*}
  as $n\rightarrow\infty$.
  \end{theorem}
  {\bf Proof:} We recall that
  \begin{eqnarray*}
      \gamma(F,t)=2(1+\overline F(t))\int_{0}^{t}\overline
      F(x)dx-2\int_{0}^{\infty}\overline F^2(x)dx+4\int_{t}^{\infty}\overline
      F^2(x)dx-\mu_F.
    \end{eqnarray*}
By the theory of Von-Mises differentiable statistical function, we have that
\begin{eqnarray}\label{vonmiss}
  \gamma(F_n,t)-\gamma(F,t)=D_1\gamma(F,t)(F_n-F)+R_n(F,t)
\end{eqnarray}
where $R_n(F,t)$ is the remainder term and $D_1\gamma(F,t)(F_n-F)$ is
G$\check{\rm a}$teaux differential of $\gamma(F,t)$ at $F$ in
direction $F_n$ and given by
\begin{eqnarray*}
  D_1\gamma(F,t)(F_n-F)=\frac{d}{d\lambda}\gamma(F+\lambda(F_n-F),t)|_{\lambda=0+}
  \end{eqnarray*}
  (see Serfling (1980), Chapter 6).
  It is easy to see that
  \begin{eqnarray*}
  D_1\gamma(F,t)(F_n-F)&=&-2(F_n(t)-F(t))\int_{0}^{t}\overline F(x)dx-2(1+\overline F(t))\int_{0}^{t}(F_n(x)-F(x))dx\\
  &&+4\int_{0}^{\infty}(F_n(x)-F(x))\overline F(x)dx-8\int_{t}^{\infty}(F_n(x)-F(x))\overline F(x)dx\\
  &&+\int_{0}^{\infty}(F_n(x)-F(x))dx.
\end{eqnarray*}
After some algebraic manipulations, we obtain
\begin{eqnarray*}
R_n(F,t)&=&\gamma(F_n,t)-\gamma(F,t)-D_1\gamma(F,t)(F_n-F)\\
&=&-2\int_{0}^{\infty}(\overline F_n(x)-\overline
F(x))^2dx+4\int_{t}^{\infty}(\overline F_n(x)-\overline F(x))^2dx\\
&&+2(\overline F_n(t)-\overline F(t))\int_{0}^{t}(\overline F_n(x)-\overline F(x))dx.
\end{eqnarray*}
For every $t\ge0$,
\begin{eqnarray*}
  \sqrt{n}|R_n(F,t)|\le 8\sqrt{n}\sup|F_n(x)-F(x)|\int_{0}^{\infty}|F_n(x)-F(x)|dx.
\end{eqnarray*}
By the classical results on the Kolmogorov-Smirnov statistics, we have
$\sqrt{n}\sup|F_n(x)-F(x)|=O_p(1)$.
Also, $\int_{0}^{\infty}|F_n(x)-F(x)|dx=o_p(1)$ (Hawkins et al. (1992), page 288). Thus $\sqrt{n}R_n(F,t)=o_p(1)$.

Let us define the stochastic process $Z_n(F)=\{Z_n(u,F): 0\le u\le 1\}$ where
\begin{eqnarray*}
 Z_n(u;F)=n^{\frac{1}2}D_1\gamma(F,F^{-1}(u))(F_n-F).
\end{eqnarray*}
By (\ref{vonmiss}), under the null hypothesis,
\begin{eqnarray}
 \sqrt{n} \gamma(F_n,F^{-1}(u))=Z_n(u;F)+\sqrt{n}R_n(F,F^{-1}(u)), \ \ \ 0\le u\le1,
\end{eqnarray}
where $F^{-1}$ is the quantile function corresponding to $F$.
Thus, to obtain the required result, it suffices to show
that under the null hypothesis
\begin{eqnarray*}
  \overline X^{-1}Z_n(u;F)\stackrel D\rightarrow Z(u), \ \ \ \ 0\le u\le 1
\end{eqnarray*}
where
$Z(u)$ is a mean-zero Gaussian process with covariance
(\ref{Covariance}). Now, we have that
\begin{eqnarray*}
n^{-\frac{1}2}Z_n(u;F)&=&-2(F_n(F^{-1}(u))-u)\int_{0}^{F^{-1}(u)}\overline F(x)dx-2(2-u)\int_{0}^{F^{-1}(u)}(F_n(x)-F(x))dx\\
  &&+4\int_{0}^{\infty}(F_n(x)-F(x))\overline F(x)dx-8\int_{F^{-1}(u)}^{\infty}(F_n(x)-F(x))\overline F(x)dx\\
  &&+\int_{0}^{\infty}(F_n(x)-F(x))du.
\end{eqnarray*}
Let us define $W_n(x)=n^{\frac12}(F_n(F^{-1}(x))-x)$, $0\le x\le 1$.
Now, under the null hypothesis $H_0:$ F is an exponential distribution
with mean $\beta$,
\begin{eqnarray*}
\beta^{-1}Z_n(u;F)&=&-2uW_n(u)-2(2-u)\int_{0}^{u}\frac{W_n(x)}{1-x}dx+4\int_{0}^{1}W_n(x)dx\\
&&-8\int_{u}^{1}W_n(x)dx+ \int_{0}^{1}\frac{W_n(x)}{1-x}dx
\end{eqnarray*}
Since $W_n(x),$ $0\le x\le1$ converges to a Browanian
bridge, $W(x)$, with covariance $E(W(s)W(t))=s(1-t)$, $s\le t$
(see Serfling (1980), p. 110) and $\overline X\stackrel {a.s.}
\rightarrow \beta$, the required result follows. \qed


Now, in order to find the critical points of $\gamma^*(F_n)$ and
\begin{eqnarray}\label{Critical}
  P\{Z>c\}=P\{\mathop{\sup}_{0\le u\le1}Z(u)>c\},
\end{eqnarray}
we use Durbin's (1985) approximation. Suppose that $y(u)$ is a
continuous Gaussian process with covariance function $\rho(s,t)$,
$0\le s\le t$ and let
\begin{eqnarray*}
 T_{u_0}=\mathop {\inf}_{u>u_0}\{u: y(u)>c(u)\}
\end{eqnarray*}
denote the first-passage time of $y(u)$ to a boundary $c(u)$ at
time $u=u_0$. Durbin (1985) approximated the density function of
the first-passage, $g(t)$, by
\begin{eqnarray}
  g(t)\simeq b(t)f(t)
\end{eqnarray}
where
\begin{eqnarray}
  b(t)=\frac{c(t)}{\rho(t,t)}\frac{\partial\rho(s,t)}{\partial
  s}|_{s=t}-c'(t).
\end{eqnarray}
and
\begin{eqnarray*}
  f(t)=(2\pi\rho(t,t))^{\frac{-1}{2}}\exp\{\frac{-c^2(t)}{2\rho(t,t)}\}.
\end{eqnarray*}

Now, suppose $c(u)=c$, then using Durbin's approximation we can
approximate (\ref{Critical}) for large value of $c$, the
probability that $Z(u)$ crosses $c$ in $[0, 1]$, as follows:

\begin{eqnarray*}
  P\{Z>c\}&=&P\{\mathop{\sup}_{0\le u\le1}Z(u)>c\}\\
  &=&P(0\le T_0\le 1)\\
  &=&\int_{0}^{1}g(t)dt\nonumber
\end{eqnarray*}
where
\begin{eqnarray*}
  g(t)\simeq \frac{6\sqrt{3}ct^2}{\sqrt{2\pi}}\exp\{-\frac{3c^2}{2}\}
\end{eqnarray*}
and thus
\begin{eqnarray}\label{Critical1}
 P\{Z>c\}\simeq\frac{2\sqrt{3}c}{\sqrt{2\pi}}\exp\{\frac{-3c^2}{2}\}.
\end{eqnarray}

In Table \ref{CriticalT}, using (\ref{Critical1}), we provide some approximated values
of critical points of $\gamma^*(F_n)$, $c$, such that $P\{Z>c\}=\alpha$, for some commonly used levels of significance.
\begin{table}[] \caption{ Approximated critical points of $\gamma^*(F_n)$ }\label{CriticalT}
 \centering
 \vspace{0.5cm}
\begin{tabular}{|l|lcccc|}

\hline
 $\alpha$&&0.01&0.025&0.05&0.1\\
 \hline
$c$&&1.9298&1.7453&1.5878&1.4065\\
 \hline
\end{tabular}
\end{table}

As for the problem of testing $H_0$ against $H'_1$ given as in (\ref{DIMRLH}),  the measure of departure from $H_0$ in favor of $H'_1$ is 
\begin{eqnarray}
\kappa(F)=\sup\{-\gamma(F,t); 0\le t<\infty\}.
\end{eqnarray}
 Similar to the test for $H_0$ against $H_1$,  we propose to reject $H_0$ in favor of $H'_1$ for large values of
\begin{eqnarray*}
\kappa^*(F_n)=\frac{\sqrt{n}\max_{0\le k\le
  n}^{}\eta_{n,k}}{\overline{X}}
\end{eqnarray*}
where
\begin{eqnarray*}
\eta_{n,k}=\sum_{i=0}^{n-1}C _{k,i}D_{i}+\sum_{i=k}^{n-1}E_{k,i}D_i
\end{eqnarray*}
and
\begin{eqnarray*}
C_{k,i}=2(1-\frac in)^2+(1-\frac in) - 2(2-\frac kn)(1-\frac in), \ E_{k,i}=2(2-\frac kn)(1-\frac in)-4(1-\frac in)^2. 
\end{eqnarray*}
From the proof of Theorem \ref{Th1}, it can be seen that $\gamma^*(F_n)$ and $\kappa^*(F_n)$ have identical asymptotic distribution. Thus, the  values given in Table \ref{CriticalT} can be used as the critical points of the test statistic $\kappa^*$ for large sample sizes.


\section{A simulation study }
In this section, we study the performance of the proposed tests. First, to investigate the speed of convergence of the test statistics $\gamma^*(F_n)$ and $\kappa^*(F_n)$ to $Z$ and accuracy the approximated critical points given in Table \ref{CriticalT}, we obtained the empirical sizes at some different nominal sizes (levels of significance). Since our test statistics are scale invariant, we generated 10000 random samples from exponential distribution with mean 1. The values, presented in Table \ref{Table-Size}, are the fraction of times that $H_0$ is rejected at some levels of significance, i.e., the test statistic is greater than the asymptotic critical value at the given nominal sizes. The results indicate that the tests have type I error rates far from the nominal sizes for small sample sizes. So, in Table \ref{Table-quantile} below, we obtained the critical points of $\gamma^*$ and $\kappa^*$ by simulation for small sample sizes $n=10(5)70$ and some commonly used levels of significance. To compute the critical points, 10000 samples were generated from the exponential distribution.     
\begin{table}[H] \caption{Empirical sizes of the tests  $\gamma^*(F_n)$ and $\kappa^*(F_n)$.}\label{Table-Size}
 \hrule
\centering 
{\fontsize{10}{9}\selectfont
\begin{tabular}{ccccccccccccc}
&&&$\gamma^*$&&&&&$\kappa^*$&\\
\cline{3-6}\cline{8-11}
 $n$ &$\alpha$&$0.01$& $0.025$& $0.05$&$0.1$&&$0.01$& $0.025$& $0.05$&$0.1$\\
\hline\\
$10$          &&0.0089    &0.0292    &0.0654    &0.1447  &&0.0013    &0.0050    &0.0119    &0.0341\\   
$15$          &&0.0107    &0.0323    &0.0681    &0.1478  &&0.0020    &0.0055    &0.0150    &0.0407\\                     
$20$          &&0.0124    &0.0346    &0.0759    &0.1512  &&0.0024    &0.0080    &0.0197    &0.0495\\                       
$25$          &&0.0113    &0.0310    &0.0665    &0.1409  &&0.0021    &0.0095    &0.0235    &0.0540\\        
$30$          &&0.0139    &0.0366    &0.0699    &0.1450  &&0.0036    &0.0101    &0.0224    &0.0565\\ 
$35$          &&0.0123    &0.0332    &0.0719    &0.1415  &&0.0042    &0.0122    &0.0260    &0.0621\\               
$40$          &&0.0134    &0.0349    &0.0719    &0.1429  &&0.0044    &0.0119    &0.0288    &0.0649\\ 
$45$          &&0.0129    &0.0336    &0.0693    &0.1403  &&0.0040    &0.0119    &0.0268    &0.0654\\ 
$50$          &&0.0138    &0.0375    &0.0686    &0.1405  &&0.0052    &0.0123    &0.0291    &0.0673\\ 
$55$          &&0.0141    &0.0338    &0.0700    &0.1392  &&0.0051    &0.0136    &0.0302    &0.0690\\ 
$60$          &&0.0124    &0.0310    &0.0661    &0.1352  &&0.0048    &0.0132    &0.0305    &0.0737\\ 
$65$          &&0.0142    &0.0325    &0.0674    &0.1379  &&0.0045    &0.0142    &0.0310    &0.0711\\ 
$70$          &&0.0127    &0.0350    &0.0706    &0.1443  &&0.0050    &0.0143    &0.0322    &0.0717\\
$100$        &&0.0122    &0.0336    &0.0649    &0.1330  &&0.0060    &0.0168    &0.0377    &0.0830\\
$200$        &&0.0115    &0.0287    &0.0645    &0.1233  &&0.0066    &0.0201    &0.0423    &0.0879\\
\hline            
\end{tabular}
}
\end{table}

\begin{table}[H] \caption{$(1-\alpha)-$quantile of $\gamma^*$ and $\kappa^*$.}\label{Table-quantile}
  \hrule
  \hrule
\centering {\fontsize{9}{9}\selectfont
\begin{tabular}{ccccccccccccc}
&&&$\gamma^*$&&&&&$\kappa^*$&\\
\cline{3-6}\cline{8-11}
 $n$ &$\alpha$&$0.01$& $0.025$& $0.05$&$0.1$&&$0.01$& $0.025$& $0.05$&$0.1$\\                                                                           
\hline\\
$10$&         &1.9179    &1.7747    &1.6429    &1.5005  &&1.5908    &1.3947    &1.2277    &1.0702\\    
$15$&        &1.9389    &1.8013    &1.6584    &1.5039  &&1.5922    &1.4050    &1.2461    &1.0846\\                                      
$20$&         &1.9636    &1.8157    &1.6755    &1.5206  && 1.5839    &1.3979    &1.2532    &1.0792\\                       
$25$&         &1.9501    &1.7870    &1.6485    &1.4995  &&1.6283    &1.4172    &1.2823    &1.1028\\        
$30$&         &1.9759    &1.8343    &1.6587    &1.5014  &&1.5976    &1.4318    &1.2830    &1.1102\\ 
$35$&         &1.9659    &1.8015    &1.6617    &1.5064  &&1.5671    &1.4247    &1.2971    &1.1242\\ 
$40$&         &1.9759    &1.8208    &1.6655    &1.5066  &&1.6197    &1.4384    &1.3033    &1.1402\\   
$45$&         &1.9758    &1.8103    &1.6656    &1.4897  &&1.6513    &1.4612    &1.3184    &1.1605\\ 
$50$&         &1.9825    &1.8180    &1.6734    &1.5000  &&1.6480    &1.4697    &1.3204    &1.1581\\ 
$55$&         &1.9720    &1.7996    &1.6669    &1.4961  &&1.6871    &1.4942    &1.3496    &1.1894\\ 
$60$&         &1.9672    &1.7856    &1.6495    &1.4855  &&1.6722    &1.4921    &1.3504    &1.1840\\ 
$65$&         &1.9834    &1.7971    &1.6541    &1.4943  &&1.6507    &1.4942    &1.3654    &1.2063\\ 
$70$&         &1.9719    &1.8056    &1.6675    &1.5021  &&1.6878    &1.5247    &1.3740    &1.2044\\
\hline                    
\end{tabular}
}
\end{table}

Since there is no other test for our considered problem and the IDMTTF (DIMTTF) test may also be used to test exponentiality against BFR (UBFR) alternative, it would be beneficial to compare it  to some  BFR (UBFR) tests. Na et al. (2005) proposed a test for the problem of testing exponentiality against BFR property based on the measure 
\begin{eqnarray}\label{NaST}
T(F)=\sup\{\phi(x; F): x\ge0\}
\end{eqnarray}
where 
\begin{eqnarray*}
\phi(x; F)=(1-\overline{F}(x))\int_{0}^{x}\overline{F}(t)dt-2\int_{0}^{x}\overline{F}^2(t)dt-\overline{F}(x) \int_{x}^{\infty}\overline{F}(t)dt + 2\int_{x}^{\infty}\overline{F}^2(t)dt.
\end{eqnarray*}
Their test statistic is 
\begin{eqnarray}
T^*=\frac{\sqrt{n}T(F_n)}{\overline{X}} 
\end{eqnarray}
where $T(F_n)$ is obtained by replacing $F$ by $F_n$ in (\ref{NaST}). There is a slight mistake in the formula for $T(F_n)$ in Na et al. (2005).  It is not difficult to see that 
\begin{eqnarray*}
T(F_n)=\max_{0\le k\le n}^{} \{\eta_1(k)-2\eta_2(k)+\eta_3(k)\}
\end{eqnarray*}
where for $k=0, \ldots, n$,
\begin{eqnarray*}
\eta_1(k)&=&\sum_{i=1}^{k}\frac{n-i+1}{n}D_{i-1},\\
\eta_2(k)&=&\sum_{i=k+1}^n\left[(\frac{n-k}{n})(\frac{n-i+1}{n})-2(\frac{n-i+1}{n})^2\right]D_{i-1},\\
\eta_3(k)&=&\sum_{i=1}^n\left[(\frac{n-k}{n})(\frac{n-i+1}{n})-2(\frac{n-i+1}{n})^2\right]D_{i-1}.
\end{eqnarray*}   
The large values of $T^*$ reject the exponentiality against BFR property. Na et al. (2005) have also proposed to reject exponentiality against UBFR property for large values of 
\begin{eqnarray}
U^*=\frac{\sqrt{n}\max_{0\le k\le n}^{} \{2\eta_2(k)-\eta_1(k)-\eta_3(k)\}}{\overline{X}}. 
\end{eqnarray}

Another well known BFR test is the one proposed by  Aarset (1985) which rejects $H_0$ in favor of BFR property for large values of
\begin{eqnarray*}
G_n=V_n+n-M_n
\end{eqnarray*}
where 
\begin{eqnarray*}
V_n&=&\min\{1\le i\le n: U_i\ge\frac in\},\\
M_n&=&\max\{0\le i\le n-1: U_i\le\frac in\}
\end{eqnarray*} 
and 
\begin{eqnarray}
U_i=\frac{\sum_{j=1}^{i}(n-j+1)D_{j}}{\sum_{j=1}^{n}(n-j+1)D_{j}}, \ \ i=1, \ldots, n, \ \ U_0=0.
\end{eqnarray}
The small values of $G_n$ reject exponentiality in favor of UBFR.
Aarset (1985) has derived the exact null distribution of $G_n$ which takes integer values in $[2, n+1]$.

To compare three tests, we use the simulated power of the tests. For a fair comparison,  we use the simulated critical points of our test and the test proposed by Na et al. (2005). We also employ the randomized version of the test of Aarset (1985).  In the alternative hypothesis, we first consider the exponential power
distribution which has been studied by Smith and Bain (1975),
Dhillon (1981), Paranjpe and Rarjarshi (1986).
The exponential power model has survival function and failure rate
function given by
\begin{eqnarray*}
  \overline F(x)=\exp\{-(e^{(\lambda x)^\beta}-1)\}
\end{eqnarray*}
and
\begin{eqnarray*}
  r(t)=\lambda\beta (\lambda t)^{\beta-1}e^{(\lambda t)^\beta},
\end{eqnarray*}
respectively. For $\beta<1$, $r(t)$  yields a bathtub shape ( Lai and Xie, 2006). Since BFR implies IDMTTF (Izadi et al., 2018), thus the exponential power is IDMTTF for $\beta<1$.
 In Table \ref{Table-PowerIDMTTF}, we simulate the power of the IDMTTF test ($\gamma^*$), Na et al.'s test ($T^*$) and Aarset's test ($G_n$) at significance levels $\alpha=0.01, 0.05$ and  $0.1$ by generating 10000 samples with sizes $n=10(10)60$ and  from the exponential power distribution with $\lambda=1$ and some values of the parameter $\beta<1$. 
The results indicate that our test and the test of Na et al. (2005) dominate Aarset's test for all cases in the alternative. For the values of $\beta$ closed to $1$, the test of  Na et al. (2005) performs better than our test,  while our test is more powerful for the other cases.

The random variable $X$ has a lognormal distribution with parameters $\mu$ and $\sigma$ if $\ln X$ has a normal distribution with mean $\mu$ and variance $\sigma^2$.  It is known that the lognormal distribution is UBFR which is also DIMTTF (Izadi et al., 2018).  In order to compare our DIMTTF test ($\kappa^*$)  with respect to  Na et al.'s test ($U^*$) and Aarset's test ($G_n$), we consider the lognormal distribution in the alternative.  In Table \ref{Table-PowerDIMTTF} below, we simulate the power of the tests at significance levels $\alpha=0.01, 0.05$ and $0.1$ by generating 10000 samples with sizes $n=10(10)60$   from the lognormal distribution with $\mu=0$ and some selected values of $\sigma$. From Table \ref{Table-PowerDIMTTF}, a similar observation to that reported for Table \ref{Table-PowerIDMTTF} can be made. The performance of Aarset test is poor. The new proposed test  performs better than the test due to Na et al. for some values of $\sigma$ and the result is reversed for other values. 

\begin{sidewaystable}[http] \caption{ The simulated power of $\gamma^*$, $G_n$, $T^*$ for the exponential power model in the alternative.}\label{Table-PowerIDMTTF}
\begin{center}
{\fontsize{9}{11}\selectfont
\begin{tabular}{|cc|ccc|ccc|ccc|ccc|ccc|}
\cline{1-17}
                  &  &$\beta$=0.1&&& $\beta$=0.3&&& $\beta$=0.5&&&$\beta$=0.7&&& $\beta=0.9$&&\\
                   \cline{3-17}
n& $\alpha$&$\gamma^*$&$G_n$& $T^*$& $\gamma^*$&$G_n$& $T^*$& $\gamma^*$&$G_n$&$T^*$ &$\gamma^*$&$G_n$&$T^*$&$\gamma^*$&$G_n$& $T^*$\\               
\cline{1-17}
10        &0.01&       0.9805 &0.0232 &0.9334&0.5118   & 0.0215& 0.2540&0.0815&0.0161&0.0143& 0.0118&0.0118 &  0.0090 &0.0229&0.0128 & 0.0497\\                        
            &0.05&      0.9950 &0.1161 &0.9814&0.7129  &0.1075 &0.4767&0.2230&0.0803& 0.0637&0.0592&0.0592& 0.0528 &0.1015&0.0641&0.1647 \\
            &0.1&      0.9976 &0.2322 &0.9908& 0.8019 &0.2150 & 0.6000&0.3232&0.1605& 0.1290&0.1183 &0.1183 & 0.1064&0.1825&0.1281 &0.2671\\
20         &0.01&      1.000&0.0415&0.9999&0.8503&0.0386& 0.7738&0.1581&0.0253 &0.0794&0.0291&0.0144 & 0.0117 &0.0523&0.0195&0.0838 \\   
             &0.05&     1.000 &0.2073&1.000&0.9392&0.1929& 0.9099 &0.3437&0.1267&0.2290&0.1063&0.0718 &0.0643 &0.1743&0.0973 &0.2697\\
             &0.1&     1.000 &0.4147&1.000&0.9656&0.3859&0.9507&0.4639&0.2534& 0.3425&0.1867&0.1437 &  0.1415&0.2788&0.1947& 0.4030\\ 
30         &0.01&      1.000&0.0598&1.000&0.9642&0.0559& 0.9515&0.2535&0.0352& 0.1737&0.0355&0.0165 & 0.0138 &0.0845&0.0281& 0.1401\\                    
             &0.05&      1.000&0.2989&1.000&0.9916&0.2797& 0.9889&0.4783&0.1758&0.3878&0.1382&0.0824 &0.0765 &0.2611&0.1406& 0.3687\\
             &0.1&        1.000&0.5977&1.000&0.9954&0.5595&0.9939 &0.6113&0.3516&0.5235&0.2385&0.1649 &0.1556 &0.3911&0.2811&  0.5177\\
40          &0.01&      1.000&0.0781&1.000&0.9929&0.0735&0.9935&0.3514&0.0452&0.3053&0.0453&0.0176 &0.0194 &0.1340&0.0373 &0.2185\\    
              &0.05&      1.000&0.3905&1.000&0.9990&0.3673&0.9990&0.5831&0.2260&0.5353&0.1527&0.0880 & 0.0918&0.3399&0.1864 &0.4632 \\
              &0.1&      1.000&0.7811&1.000&0.9994&0.7345&0.9994&0.7037&0.4521&0.6538&0.2556&0.1760 &0.1788 &0.4865&0.3729 &0.5926\\         
50          &0.01&      1.000&0.0965&1.000&0.9983&0.0904&0.9982&0.4373&0.0532&0.4070&0.0536&0.0195 & 0.0181&0.1727&0.0452 &0.2764\\               
              &0.05&      1.000&0.4824&1.000&1.0000&0.4522& 0.9999&0.6657&0.2661&0.6561&0.1838&0.0974 &0.1002 &0.4133&0.2261 &0.5443\\        
              &0.1&      1.000&0.9649&1.000&1.0000&0.9044&1.0000&0.7821&0.5322&0.7671&0.3037&0.1948 &0.2039 &0.5777&0.4522 & 0.6870\\
60          &0.01&      1.000&0.1149&1.000&1.000&0.1075&1.000&0.5410&0.0632&0.5167&0.0716&0.0224 & 0.0194& 0.2448&0.0543 & 0.3480\\               
              &0.05&      1.000&0.5745&1.000&1.000&0.5376&1.000&0.7623&0.3159&0.7424&0.2232&0.1119 &0.1093 &0.5129&0.2716 &0.6233\\        
              &0.1&        1.000&1.000&1.000&1.000&0.9874&1.000&0.8523&0.6350&0.8373&0.3517&0.2276 & 0.2212&0.6651&0.5458
 &0.7529\\        
 \cline{1-17}
\end{tabular}}
\end{center}
\end{sidewaystable}

\begin{sidewaystable}[http] \caption{ The simulated power of $\kappa^*$, $G_n$ and $U^*$ for the log-normal distribution in the alternative.}\label{Table-PowerDIMTTF}
\begin{center}
{\fontsize{9}{11}\selectfont
\begin{tabular}{|cc|ccc|ccc|ccc|ccc|ccc|}
\cline{1-17}
                    &  &$\sigma$=0.2&&& $\sigma$=0.6&&& $\sigma$=1&&&$\sigma$=1.4&&& $\sigma=2$&&\\
                   \cline{3-17}
n& $\alpha$&$\kappa^*$&$G_n$& $U^*$& $\kappa^*$&$G_n$& $U^*$& $\kappa^*$&$G_n$&$U^*$ &$\kappa^*$&$G_n$&$U^*$&$\kappa^*$&$G_n$& $U^*$\\                 
\cline{1-17}
10        &0.01&  1&0.000& 1   &0.2860& 0.0152  & 0.2868 &0.0270&  0.0249 & 0.0202 &0.0738 &0.0100 &0.0562&0.2878& 0.0016&0.2751         \\                             
            &0.05&  1&0.000& 1   & 0.6166& 0.0760 & 0.6139 &0.1179&0.1245&0.0871&0.1886& 0.0501& 0.1455&0.4785& 0.0080&0.4558   \\
            &0.1&   1&0.000& 1   & 0.7460&0.1521  & 0.7587 &0.1912&0.2489&0.1567 &0.2619& 0.1001 &0.2100&0.5684&0.0160 &0.5484      \\
20         &0.01 &1&0.000&1&0.6461&0.0176&0.6800& 0.0768 &  0.0345&0.0371&0.3139 &0.0135 &0.2644&0.7610 &   0.0006&0.7574 \\                   
             &0.05&1&0.000&1&0.8749&0.0878&0.9031& 0.2063&  0.1723&0.1328&0.4748 &0.0677 & 0.4187&0.8645& 0.0031&0.8668\\
             &0.1 &1&0.000&1&0.9428&0.1756&0.9571 &0.3204 & 0.3446& 0.2349 &  0.5702&0.1353 & 0.5120 &0.9069 & 0.0062& 0.9053\\ 
30         &0.01&1&0.000&1&0.8620& 0.0173 &0.8884& 0.1406&0.0386&0.0613&0.5044&0.0173&0.4543& 0.9311&0.0004& 0.9339   \\                    
             &0.05&1& 0.000&1&0.9655& 0.0867&0.9765& 0.3038&0.1932&0.1898&0.6602&0.0864&0.6214&0.9680 &0.0020& 0.9712     \\
             &0.1& 1&0.000 &1&0.9878&0.1734& 0.9918& 0.4245&0.3864&0.3084&0.7333&0.1727&0.7007&0.9792& 0.0039&0.9803    \\
40          &0.01&1& 0.000 &1&0.9507&0.0186&0.9619&0.1946&0.0419& 0.0809&0.6585&0.0207&0.6120&0.9836&0.0003&0.9841\\    
              &0.05&1&0.000&1&0.9934&0.0931&0.9960&0.3883&0.2094&0.2398&0.7868& 0.1036& 0.7650&0.9939&0.0016 &0.9943 \\
              &0.1&1& 0.000&1&0.9977&0.1861& 0.9992&0.5153&0.4188&0.3776& 0.8486&0.2071&0.8307&0.9964 & 0.0032&0.9969 \\         
50          &0.01&1&0.000&1&0.9853& 0.0186 &0.9918&0.2548&0.0437&0.1205 &0.7568&0.0236&0.7368& 0.9961& 0.0002&0.9974 \\               
              &0.05&1&0.000&1&0.9991&0.0930 &0.9994&0.4607&0.2186& 0.2898& 0.8660&0.1180&0.8538& 0.9991&0.0011&0.9991    \\        
              &0.1&1&0.000&1&0.9997& 0.1860 &0.9998& 0.5803&0.4372& 0.4261 &0.9077&0.2361&0.8977&0.9995&0.0021&0.9997  \\
60          &0.01&1&0.000&1&0.9958&0.0191&0.9983&0.3018&0.0445&0.1624&0.8338&0.0265&0.8412& 0.9988&0.0000&0.9991         \\               
              &0.05&1&0.000&1&0.9996&0.0957&0.9999& 0.5135&0.2227&0.3645&0.9171&0.1324&0.9196&0.9997&0.0004&0.9998      \\        
              &0.1& 1&0.000&1&1.0000&0.1914&1.0000&0.6407&0.4454&0.5191&0.9446&0.2649& 0.9485&0.9999& 0.0008&1.0000        \\        
 \cline{1-17}
\end{tabular}}
\end{center}
\end{sidewaystable}


\newpage
  \center {\bf References}
  \small{
\begin{description}
\item Aarset, M.V., 1985. The null distribution for a test of constant versus ``bathtub" failure rate. \emph{Scand. J. Statist.} { 12}, 55-61.

\item Anis, M.Z., 2014. Tests of non-monotonic stochastic aging notions in reliability theory. \emph{Stat. Papers}.   55, 691--714.

\item Anis, M.Z., Ghosh, A., 2015. Monte Carlo comparison of tests of exponentiality against NWBUE alternatives. \emph{Math. Comput. Simulation}. 115, pp.1-11. 


\item Asha, G., Nair, N.U., 2010. Reliability properties of mean time to failure in age replacement models. \emph{Int. J.
Reliab. Qual. Saf. Eng.} {17}, 15--26.



\item Barlow, R.E., Proschan, F., 1965. Mathematical theory of reliability. Wiley, New York


\item Belzunce, F., Ortega, E. and Ruiz, J.M., 2007. On non-monotonic ageing properties from the Laplace transform, with
actuarial applications. \emph{Insurance Math. Econom.} {40}, 1--14.




\item Deshpande,  J.V.,   Suresh,   R.P.,  1990.  Non-monotonic  ageing.  \emph{ Scand.  J.  Statist.}  17,  257--262. 

\item Dhillon, B., 1981. Life distributions. \emph{IEEE Trans. Reliab.} { R-30}, 457-460.

\item Durbin, J., 1985. The first-passage density of a continuous Gaussian process to a general boundary. \emph{ J. Appl. Probab.}
22, 99--122.

\item Glaser, R.E., 1980. Bathtub and related failure rate characterizations. \emph{J. Amer. Statist. Assoc.} 75, 667--672.

\item Guess, F.,  Hollander, M.,  Proschan, F., 1986. Testing exponentiality versus a trend change in mean residual life. \emph{Ann. Statist.} {14}. 1388-1398.

\item Hawkins, D.L., Kochar, S., 1997. Inference about the transition-point in NBUE-NNWUE or NWUE-NBUE. \emph{Sankhy\={a} Series A.} 59, 117-132.

\item  Hawkins, D.L., Kochar, S., Loader, C., 1992. Testing exponentiality against IDMRL distributions with unknown change point, \emph{Ann. Statist.} 20, 280--290.

\item Izadi, M., Sharafi, M., Khaledi, B.H., 2018. New non-parametric classes of distributions in terms of mean time to failure in age replacement. \emph{J. Appl. Probab.} accepted for publication.

\item Kattumannil, S.K., Anisha, P., 2016. A simple non-parametric test for decreasing mean time to failure. \emph{Stat. Papers.} doi: 10.1007/s00362-016-0827-y.

\item Kayid, M., Ahmad, I.A., Izadkhah S., Abouammoh, A.M., 2013. Further results involving the mean time to
failure order, and the decreasing mean time to failure class. \emph{IEEE Trans. Reliab.} {62}, 670--678.



\item Klefsj\"{o}, B.,  1982.  On  aging  properties  and  total  time  on test  transforms. 
\emph{Scand. J. Statist.}  {9},  37--41. 

\item Klefsj\"{o}, B., 1989. Testing against a change in the NBUE property. \emph{Microelectron. Reliab.} {29}, 559-570.

\item Knopik, L., 2005. Some results on the ageing class. \emph{ Control Cybern}. {34}, 1175--1180.

\item Knopik, L., 2006. Characterization of a class of lifetime distributions.  \emph{Control Cybern.} { 35}, 407--414.


\item Lai, C.D., 1994. Tests of univariate and bivariate stochastic ageing. \emph{IEEE Trans. Reliab.} { R-43}, 233--241.

\item Lai, C.D., Xie, M., 2006. Stochastic ageing and dependence for reliability. Springer, New York. 


\item Li, X.,  Xu, M., 2008. Reversed hazard rate order of equilibrium distributions and a related ageing notion. \emph{Stat. Papers} {49}, 749--767.


\item Mitra, M., Basu, S.K., 1994. On a nonparametric family of life distributions and its dual. \emph{J. Statist. Plann. Inference} 39, 385--397.

\item Na, M.H., Jeon, J. and Park, D.H., 2005. Testing whether failure rate changes its trend with unknown change points. \emph{J. Statist. Plann. Inference}, {129}, 317--325.



\item Nakagawa, T., 2005. Maintenance theory of reliability. London: Springer

\item Paranjpe, S.A.,  Rajarshi, M.B., 1986. Modelling non-monotonic survivorship data with bathtub distributions. \emph{Ecology}. {67}, 1693-1695.

\item Rajarshi, S. and Rajarshi, M. B. (1988), Bathtub distributions: a review. \emph{ Comm. Statist. Theory Methods} 17, 2597--2621.







\item Serfling, R.J., 1980. Approximation Theorems of Mathematical
Statistics. Wiley, New York.

\item Smith, R.M., Bain,  L.J., 1975. An exponential power life-testing distribution. \emph{Commun. Statist. - Theor. Meth.}  { 4}, 469-481.


\end{description}
}

\end{document}